%%%%%%%%%%%%%%%%%%%%%%%%%%%%%%%%%%%%%%%%%%%%%%%%%%%%%%%%%%%%%
%%  On a Claim of Ramanujan about Certain Hypergeometric Series
%%
%%                 author David Bradley
%%
%%           typed by Hilda K. Britt May 1992
%%%%%%%%%%%%%%%%%%%%%%%%%%%%%%%%%%%%%%%%%%%%%%%%%%%%%%%%%%%%%
\documentstyle{amsppt}
\magnification=1200
\NoBlackBoxes
%\TagsOnRight
\NoRunningHeads

\catcode`\@=11
\font@\twelvesmc=cmcsc12
\def\headfont@{\twelvesmc}
\catcode`\@=\active

\define\lrangle#1{{\langle #1 \rangle}}
\define\lrfloor#1{{\lfloor #1 \rfloor}}
\define\lrceil#1{{\lceil #1 \rceil}}

\redefine\Re{\operatorname{Re}}
\define\a{\alpha}
\redefine\b{\beta}

\define\({\left(}
\define\){\right)}
\define\[{\left]}
\define\]{\right]}

\topmatter

\title On a Claim of Ramanujan about Certain Hypergeometric
Series
\endtitle

\author David Bradley \endauthor

\abstract We state and prove a claim of Ramanujan. As a
consequence, a large new class of Saalsch\"utzian hypergeometric
series is summed in closed form.
\endabstract

\endtopmatter

\document
\baselineskip=24pt plus 2pt

Although Ramanujan published no papers on hypergeometric
series, chapters 10 and 11 of his second notebook \cite5 deal almost
exclusively with this subject. Furthermore, several other results
on hypergeometric functions are scattered among the 100 pages of
unorganized material at the end of the second notebook. In 1923,
Hardy [4] published a brief survey of chapter 10's corresponding
chapter in Ramanujan's first notebook, namely chapter 12. There
we
see that Ramanujan had rediscovered many of the classical formulae
of the subject, including the theorems of Gauss, Kummer, Dougall,
Dixon, and Saalsch\"utz. However, in addition, Ramanujan
discovered many new theorems about hypergeometric series, in
particular, theorems on products of hypergeometric series and
several types of asymptotic expansions. For proofs, see Berndt [3].
In this paper, we examine
an enigmatic claim about hypergeometric series made by Ramanujan
in the unorganized portion of the second notebook, and it is to
this claim that we now turn.

At the top of page 280 in his second notebook \cite5, Ramanujan
states ``the difference between
$$
\frac{\Gamma(\beta+1-m)}{\Gamma(\a+\b+1-m)}
$$
and
$$
\aligned
\frac{\Gamma(\b+1)}{\Gamma(a+\b+1)} +
\frac{\a m}{1!}  &\cdot
\frac{\Gamma(\b+n+1)}{\Gamma(\a+\b+n+2)} +
\frac{\a(\a+1)}{2!} \\
& \cdot
\frac{m(m+2n+1)\Gamma(\b+2n+1)}{\Gamma(\a+\b+2n+3)} +
\frac{\a(\a+1)(\a+2)}{3!} \\
& \cdot
\frac{m(m+3n+1)(m+3n+2)\Gamma(\b+3n+1)}{\Gamma(\a+\b+3n+4)}
+\cdots\text{''}
\endaligned
$$
It is not entirely clear what Ramanujan meant by this statement,
nor even what values of the parameters $\a$, $\b$, $m$, $n$ he
considered. In this paper, we give perhaps the most plausible
interpretation of Ramanujan's enigmatic entry. In fact, we
determine conditions for which this difference equals zero.

For each integer $j$ and complex $a$, $z$, $a_1,\dots,a_p$,
$b_1,\dots,b_q$, as customary, set
$$
(a)_j = \frac{\Gamma(a+j)}{\Gamma(a)}
$$
and
$$
{}_pF_q \bigg( \matrix
a_1,\dots,a_p \\
b_1,\dots,b_q  \endmatrix\bigg| z \bigg) =
\sum^\infty_{j=0}
\frac{(a_1)_j\cdots (a_p)_j z^j}
     {(b_1)_j\cdots (b_q)_j j!} \, .
$$
When $p=q+1$, ${}_pF_q$ converges at $z=1$ if
$\Re(b_1+\cdots+b_q) > \Re(a_1 + \cdots + a_p)$. See \cite{2, p.
8}. If $p=2$ and $q=1$, we have Gauss's theorem; namely,
$$
{}_2F_1\bigg( \matrix a, b \\ c \endmatrix \bigg|1 \bigg) =
\frac{\Gamma(c)\Gamma(c-a-b)}{\Gamma(c-a)\Gamma(c-b)}, \quad
\Re(c) > \Re(a+b).
\tag1
$$
See \cite{2, p. 2} for a proof.

Let $\a$, $\b$, $m$ be complex numbers. Define the function $S$,
taking values from the open complex plane into the Riemann
sphere, by
$$
S(z)=S(\a,\b,m,z) = m \sum^\infty_{j=0}
\frac{\Gamma(\b+1+jz)\Gamma(m+j(z+1))}
     {\Gamma(\a+\b+1+j(z+1))\Gamma(m+jz+1)} \cdot
\frac{(\a)_j}{j!}.
\tag2
$$
Note that $S(z)$ coincides with Ramanujan's hypergeometric series
when the variable $z$ is replaced by $n$. When the parameter $\a$
is a non--positive integer, the series terminates, and in this
case, we have the following theorem.

\proclaim{Theorem} Let $\a=-k$, where $k$ is a nonnegative
integer. Then, for every complex number $z$,
$$
S(z) = \frac{\Gamma(\b+1-m)}{\Gamma(\a+\b+1-m)}.
\tag3
$$
Furthermore, if $z=0$ and $\Re(\b+1-m)>0$, then (3) is valid for
all complex numbers $\a$.
\endproclaim

Before proving the theorem, some preliminary remarks are in
order. When $z=n$ is a nonnegative integer, the series (2)
reduces to
$$
S(n) = \frac{\Gamma(\b+1)}{\Gamma(\a+\b+1)} \sum^\infty_{j=0}
\frac{(\b+1)_{nj}(m)_{(n+1)j}(\a)_j}
     {(\a+\b+1)_{(n+1)j}(m+1)_{nj}j!},
\tag4
$$
which can be recast in the form
$$
\aligned
&S(n) = \frac{\Gamma(\b+1)}{\Gamma(\a+\b+1)}  \\
&\times {}_{2n+2}F_{2n+1} \bigg( \matrix
&\dfrac{\b+1}n, \dfrac{\b+2}n,\dots,\dfrac{\b+n}n, \dfrac m{n+1},
\dfrac{m+1}{n+1},\dots,\dfrac{m+n}{n+1}, \a \\
&\dfrac{m+1}n, \dfrac{m+2}n,\dots,\dfrac{m+n}n, \dfrac{\a+\b+1}{n+1},
\dfrac{\a+\b+2}{n+1},\dots,\dfrac{\a+\b+1+n}{n+1} \endmatrix
\bigg| 1\bigg)
\endaligned\tag5
$$
by means of Gauss's multiplication formula. Observe that
${}_{2n+2}F_{2n+1}$ in (5) is Saalsch\"utzian for $n>0$, i.e.,
the sum of the denominator parameters exceeds the sum of the
numerator parameters by $1$. Thus, for each positive integer $n$,
${}_{2n+2}F_{2n+1}$ converges for all complex $\a$, $\b$, $m$ by
the previously cited remarks on convergence.

R. Askey observed that
in the form (5), the terminating case with $n=1$ follows readily
from an unpublished result of Askey and Ismail \cite1. By
utilizing Euler's evaluation of the Beta integral and a
transformation of Pfaff, they showed that, if $\Re(d)>\Re(a)>0$
and $k$ is a nonnegative integer, then
$$
{}_4F_3 \bigg( \matrix
&\dfrac a2, \dfrac{a+1}2, -k, c \\
&\dfrac d2 , \dfrac{d+1}2, -k+a+c+1-d \endmatrix
\bigg| 1\bigg) =
\frac{(d-a)_k(d-c)_k}{(d-a-c)_k(d)_k} {}_3F_2
\bigg( \matrix &-k, a, c \\ &k+d, d-c \endmatrix \bigg| 1\bigg).
$$
As Askey observed, if $c$ tends to $k+d$, the expression on the left becomes a
$_4F_3$ of the form occurring in (5) with $n=1$. The expressions
on the right are easily evaluated using the fact that
$$
\lim_{\epsilon\to 0} {}_2F_1
\bigg( \matrix -k, a \\ -k+\epsilon \endmatrix
\bigg| 1\bigg) = \frac{(-k-a)_k}{(-k)_k}.
$$

We are grateful to the referee for pointing out that the
aforementioned  ${}_4F_3$ transformation of Askey and Ismail is a
specialization of a result due to F.~J.~W.~Whipple, and can be
obtained by replacing $k$ by $d$, $a$ by $-a+d$, $c$ by $a$, then
$b$ by $c$ and $m$ by $k$ in the formula (1) of [2, p. 32].

\demo{Proof of Theorem} The case $z=0$ is easily proved by
setting $n=0$ in (4) and applying Gauss's theorem (1).

Next, let $\a=-k$. Since
$$
\align
S(z) &= m \sum^k_{j=0}
\frac{\Gamma(\b+1+jz)\Gamma(m+j(z+1))(-k)_j}
     {\Gamma(-k+\b+1+j(z+1))\Gamma(m+jz+1)j!} \\
&= m \sum^k_{j=0} (-k+\b+1+j(z+1))_{k-j} (m+jz+1)_{j-1}(-k)_j/j!
\endalign
$$
is a polynomial in $z$ of degree $k-1$, it suffices to prove (3)
for $k$ distinct values of $z$. We shall, in fact, prove that (3)
holds for all positive integers $z$, since the argument is no
more difficult than the argument for only $k$ values of $z$. So
let $z=n$ be a positive integer. From (4), since $\a=-k$, we have
$$
\aligned
\frac{\Gamma(\a+\b+1-m)}{\Gamma(\b+1-m)} S(n)
&= \frac{\Gamma(\b+1)\Gamma(\a+\b+1-m)}{\Gamma(\a+\b+1)\Gamma(\b+1-m)}
\sum^k_{j=0}
\frac{(\b+1)_{nj}(m)_{(n+1)j}(\a)_j}
     {(\a+\b+1)_{(n+1)j}(m+1)_{nj}j!} \\
&= \frac{\Gamma(m-\b)\Gamma(-\a-\b)}{\Gamma(-\b)\Gamma(m-\a-\b)}
\sum^k_{j=0} \frac{(\b+1)_{nj}(m)_{(n+1)j}(\a)_j}
     {(\a+\b+1)_{(n+1)j}(m+1)_{nj}j!}
\endaligned\tag6
$$
by the reflection formula for $\Gamma$. The sine factors that
would normally appear reduce to unity because $\a$ is an integer.
By Gauss's theorem (1), for $\Re(-\a-\b-j(n+1))>0,$ we have
$$
\aligned
{}_2F_1 \bigg( \matrix m+j(n+1), \a+j \\ m-\b+j \endmatrix
\bigg| 1\bigg)
&= \frac{\Gamma(m-\b+j)\Gamma(-\a-\b-j(n+1))}
     {\Gamma(m-\a-\b)\Gamma(-\b-jn)} \\
&=\frac{\Gamma(m-\b)(m-\b)_j\Gamma(-\a-\b)(-\b-jn)_{nj}}
       {\Gamma(m-\a-\b)(-\a-\b-j(n+1))_{(n+1)j}\Gamma(-\b)} \\
&=\frac{\Gamma(m-\b)\Gamma(-\a-\b)}{\Gamma(-\b)\Gamma(m-\a-\b)}
\cdot \frac{(m-\b)_j(\b+1)_{nj}(-1)^j}{(\a+\b+1)_{(n+1)j}} .
\endaligned\tag7
$$
The condition $\Re(-\a-\b-j(n+1))>0$ is certainly satisfied if
$\Re(\b)<-nk$. For now, assume this. Then, (6) and (7) imply that
$$
\aligned
S(n) &= \frac{\Gamma(\b+1-m)}{\Gamma(\a+\b+1-m)} \sum^k_{j=0}
\frac{(\a)_j(m)_{(n+1)j}(-1)^j}{(m-\b)_j(m+1)_{nj}j!} {}_2F_1
\bigg( \matrix m+j(n+1), \a+j \\ m-\b+j \endmatrix\bigg| 1\bigg) \\
&=\frac{\Gamma(\b+1-m)}{\Gamma(\a+\b+1-m)}  \sum^k_{j=0}
\frac{(\a)_j(m)_{(n+1)j}(-1)^j}{(m-\b)_j(m+1)_{nj}j!}
\sum^\infty_{s=0} \frac{(m+j(n+1))_s(\a+j)_s}{(m-\b+j)_s s!} \\
&=\frac{\Gamma(\b+1-m)}{\Gamma(\a+\b+1-m)}  \sum^k_{j=0}
\frac{(-1)^j/j!}{(m+1)_{nj}} \sum^\infty_{s=0}
\frac{(m)_{s+j(n+1)}(\a)_{s+j}}{(m-\b)_{s+j} s!} \\
&= \frac{\Gamma(\b+1-m)}{\Gamma(\a+\b+1-m)}  \sum^k_{j=0}
\frac{(-1)^j/j!}{(m+1)_{nj}} \sum^\infty_{r=j}
\frac{(m)_{r+nj}(\a)_r}{(m-\b)_r(r-j)!} \\
&= \frac{\Gamma(\b+1-m)}{\Gamma(\a+\b+1-m)}  \sum^\infty_{r=0}
\frac{(\a)_r(m)_r}{(m-\b)_r r!} \sum^r_{j=0}
\frac{(m+r)_{nj}}{(m+1)_{nj}} (-1)^j \binom rj.
\endaligned\tag8
$$
Interchanging the order of summation is justified since both sums
contain only finitely many nonzero terms.

Let $E$ denote the inner sum in (8). For $r=0$, it is clear that
$E=1$. For $r\geq 1$, we shall prove that $E=0$. Note that
$$
\frac{(m+r)_{nj}}{(m+1)_{nj}} =
\frac{\Gamma(m+r+nj)\Gamma(m+1)}{\Gamma(m+1+nj)\Gamma(m+r)} =
\frac{(m+nj+1)_{r-1}}{(m+1)_{r-1}}.
$$
Thus,
$$
E=\sum^r_{j=0} \frac{(m+nj+1)_{r-1}}{(m+1)_{r-1}} (-1)^j \binom
rj.
$$
But,
$$
(k)_j = D^j x^{k+j-1}\big|_{x=1}.
$$
Therefore,
$$
\aligned
E&=\frac1{(m+1)_{r-1}} \sum^r_{j=0} D^{r-1} x^{(m+nj+1)+(r-1)-1}
(-1)^j \binom rj \bigg|_{x=1} \\
&= \frac{D^{r-1}x^{m+r-1}}{(m+1)_{r-1}} \sum^r_{j=0} x^{nj}
(-1)^j \binom rj \bigg|_{x=1} \\
&= \frac{D^{r-1}x^{m+r-1}(1-x^n)^r}{(m+1)_{r-1}}\bigg|_{x=1} \\
&=0,
\endaligned
$$
as required.

Hence, for all positive integers $n$, if $\Re(\b)<-nk$, (3)
holds. Since both sides of (3) define meromorphic functions of
$\b$ for fixed $n$, $m$ and $\a=-k$, the restriction $\Re(\b)<-
nk$ may be removed by analytic continuation. By our earlier
remarks, this completes the proof of the theorem.

In general, if $S(z)$ is nonterminating, (3) is false. For
example, let $\a$, $\b$, $m$ be any complex numbers satisfying
$\a+\b+1=m$. By (5),
$$
\aligned
S(1) &= \frac{\Gamma(\b+1)}{\Gamma(\a+\b+1)} {}_4F_3 \bigg(
\matrix &\b+1, \dfrac m2, \dfrac{m+1}2 , \a \\
        &m+1, \dfrac{\a+\b+1}2, \dfrac{\a+\b+2}2
\endmatrix \bigg| 1\bigg) \\
&= \frac{\Gamma(\b+1)}{\Gamma(m)} {}_2F_1 \bigg(
\matrix \b+1, \a \\ m+1
\endmatrix \bigg| 1\bigg) \\
&= \frac{\Gamma(\b+1)\Gamma(m+1)\Gamma(m-\b-\a)}
        {\Gamma(m)\Gamma(m-\b)\Gamma(m+1-\a)} \\
&=\frac m{(m-\a)\Gamma(\a+1)}.
\endaligned\tag9
$$
But
$$
\frac{\Gamma(\b+1-m)}{\Gamma(\a+\b+1-m)} = \frac{\Gamma(-\a)}{\Gamma(0)}
(= 0 \quad \text{for} \quad \a \ne 0,1,2,3 \dots).
\tag10
$$
When $\a$ is a non--positive integer, (9) and (10) must be equal,
as we have seen. Indeed, they both vanish when $\a$ is a negative
integer. However, it is clear that (9) and (10) are unequal for
non--integral $\a$: (10) vanishes while (9) does not. We have
used
MAPLE to calculate several other nonterminating examples and
found no instances when (3) is valid.
\enddemo

\subhead Acknowledgement \endsubhead The author would like to
thank Professor Bruce Berndt for suggesting this problem to him.

\enddocument